\documentclass{amsart}

\newtheorem{theorem}{Theorem}
\newcommand{\bt}{\begin{theorem}}
\newcommand{\et}{\end{theorem}}
\newtheorem*{theoremNN}{Theorem}
\newcommand{\btNN}{\begin{theoremNN}}
\newcommand{\etNN}{\end{theoremNN}}

\newtheorem{lemma}{Lemma}
\newcommand{\bl}{\begin{lemma}}
\newcommand{\el}{\end{lemma}}
\newtheorem{corollary}{Corollary}
\newcommand{\bc}{\begin{corollary}}
\newcommand{\ec}{\end{corollary}}
\newtheorem{example}{Example}
\newcommand{\bex}{\begin{example}}
\newcommand{\eex}{\end{example}}

\newcommand{\beq}{\begin{equation}}
\newcommand{\eeq}{\end{equation}}
\newcommand{\benum}{\begin{enumerate}}
\newcommand{\eenum}{\end{enumerate}}

\newcommand{\R}{\ensuremath{\mathbf R}}
\newcommand{\C}{\ensuremath{\mathbf C}}

\DeclareMathOperator{\colsum}{\text{colsum}}

\newcommand{\bsmallmat}{\left(\begin{smallmatrix}}
\newcommand{\esmallmat}{\end{smallmatrix}\right)}

\newcommand{\bmat}{\left(\begin{matrix}}
\newcommand{\emat}{\end{matrix}\right)}

\DeclareMathOperator{\qqand}{\qquad\text{and}\qquad}

\newtheorem{problem}{Problem}
\newcommand{\bprob}{\begin{problem}}
\newcommand{\eprob}{\end{problem}}

\numberwithin{equation}{section}

\usepackage{amssymb,latexsym}

\title[Multiplicity matrices]{Interactions of zeros of polynomials and multiplicity matrices}

\author{Melvyn B. Nathanson}

\address{Department of Mathematics\\Lehman College (CUNY)\\Bronx, NY 10468}

\email{melvyn.nathanson@lehman.cuny.edu}

\date{\today}



\thanks{Supported in part by PSC-CUNY Research Award Program grant 63117-00 51.}

\begin{document}

\begin{abstract} 
An  $m \times (n+1)$ multiplicity matrix is a matrix $M = \bmat  \mu_{i,j}  \emat$ 
with rows  enumerated by $i \in \{ 1,\ 2, \ldots, m \}$ and columns  
enumerated by $j \in \{ 0,1,\ldots, n \}$ whose coordinates are nonnegative integers 
 satisfying the following  two properties:
(1) If $\mu_{i,j} \geq 1$, then $j \leq n-1$ and $\mu_{i,j+1} = \mu_{i,j}-1$, and 
(2) $\colsum_j(M) = \sum_{i=1}^{m} \mu_{i, j} \leq n-j$ for all $j$.  

Let $K$ be a field of characteristic 0 and let $f(x)$ be a polynomial of degree $n$ 
with coefficients in $K$.  Let $f^{(j)}(x)$ be the $j$th 
derivative of $f(x)$.
Let $\Lambda = ( \lambda_1,\ldots, \lambda_{m})$ 
be a sequence of distinct elements of $K$.  
For $i \in \{1, 2, \ldots, m \}$ and  $j \in \{1,2,\ldots, n\}$, 
let $ \mu_{i,j}$ be the multiplicity of $\lambda_i$ as a zero of the polynomial $f^{(j)}(x)$. 
The $m \times (n+1)$ matrix $M_f(\Lambda) = \bmat  \mu_{i,j}  \emat$ is called the 
 \emph{multiplicity matrix of the polynomial $f(x)$ with respect to $\Lambda$}.  Conditions for a multiplicity matrix to be the multiplicity matrix of a polynomial 
 are established, and examples are constructed of multiplicity matrices that are not multiplicity matrices of polynomials.
An open problem is to classify the multiplicity matrices 
that are multiplicity matrices 
 of polynomials in $K[x]$ and to construct multiplicity matrices that are not multiplicity matrices of polynomials.   
 \end{abstract}

\maketitle

\section{Interactions of zeros}       \label{multiplicity:section:intro}
Let $K$ be a field of characteristic 0.  In this paper we consider polynomials $f(x)$ 
in the polynomial ring $K[x]$.  
   The element $\lambda \in K$ is a \emph{zero}\index{zero}  
of the  polynomial $f(x)$ 
of \emph{multiplicity}\index{multiplicity} $\mu \geq 1$ if there is a polynomial $g(x) \in K[x]$ such that 
\[
f(x) = (x-\lambda)^{\mu} g(x) \qqand g(\lambda) \neq 0.
\]
If $f(\lambda) \neq 0$, then $\lambda$ is a zero of $f(x)$ of multiplicity 0.
The multiplicity of a zero is also called the \emph{order}\index{order} of the zero.  

Let $f(x)$ be a polynomial of degree $n$ and let $f^{(j)}(x)$ be the $j$th 
derivative of $f(x)$.
Let $\Lambda = (\lambda_1, \lambda_2,\ldots, \lambda_{m})$ 
be a sequence of distinct elements of $K$.  
For $i \in \{1,2  \ldots, m \}$ and  $j \in \{1,2,\ldots, n\}$, 
let $ \mu_{i,j}$ be the multiplicity of $\lambda_i$ as a  zero of the polynomial $f^{(j)}(x)$. 
We are interested in the interactions of zeros of polynomials.  
Here is a typical question.   
For distinct  elements $\lambda_1$ and $\lambda_2$, we ask: 
To what extent do the multiplicities of the  element $\lambda_1$ as a zero 
of the polynomials in the sequence of derivatives $D_f(x) = \left( f(x), f'(x), f''(x), \ldots, f^{(n)}(x) \right)$ 
constrain the multiplicities of the  element $\lambda_2$ as a zero 
of the polynomials in the sequence $D_f(x)$?  

Consider, for example, the quartic polynomials 
\[
f_1(x) = x^4 + x^3 
\]
and
\[
f_2(x) = x^4 - 4x^3   + 7x^2 -6x+3.      
\]
Let $\Lambda = ( \lambda_1, \lambda_2 )$, where $\lambda_1 = 0$ and $\lambda_2 = 1$. 
For $j \in \{0,1,2,3,4\}$, let $\mu_{1,j}$ be the multiplicity of $\lambda_1$ 
as a zero of $f^{(j)}_1(x)$ 
and let  $\mu_{2,j}$ be the multiplicity of $\lambda_2$ 
as a zero of $f^{(j)}_2(x)$.  
We have 
\[
M_{f_1}(0) = M_{f_1}(\lambda_1) = 
( \mu_{1,0 }  ,   \mu_{1, 1}  ,  \mu_{1, 2} ,  \mu_{1,3 }  ,  \mu_{1,4 }\ ) = ( 3,2,1,0,0)
\]
and 
 \[
M_{f_2}(1) = M_{f_1}(\lambda_2) = 
(\mu_{2, 0}, \mu_{2,1 }, \mu_{2, 2}, \mu_{2, 3}, \mu_{2, 4}) = ( 0,1,0,1,0).  
\]
In Section~\ref{multiplicity:section:MM} we prove that there exists no quartic 
polynomial $f(x)$ and no sequence $\Lambda = (\lambda_1, \lambda_2)$ such that 
the multiplicity of $\lambda_i$ as a zero of $f^{(j)}(x)$ is $\mu_{i,j}$ 
for all $i \in \{1,2\}$ and $j \in \{0,1,2,3,4\}$.  
However, there does exist a unique monic quintic polynomial 
\[
g(x) = x^5 -\frac{25}{8}x^4 + \frac{5}{2} x^3 
= (x-1)^5  + \frac{15}{8} (x-1)^4 - \frac{5}{4} (x-1)^2 + \frac{3}{8}
\]
that satisfies the requirement that the multiplicity of $\lambda_i$ as a zero of $g^{(j)}(x)$ is $\mu_{i,j}$ 
for all $i \in \{1,2\}$ and $j \in \{0,1,2,3,4\}$. 
This example is related to questions~\eqref{Multiplicity:problem:Np1} 
and~\eqref{Multiplicity:problem:Np2} in Section~\ref{multiplicity:section:OpenProblems}.

The multiplicity problems in this paper arose from a study of the  Budan-Fourier theorem.  

\btNN[Budan-Fourier] 
Let $f(x)$ be a polynomial of degree $n$ with real coefficients, and let 
\[
D_f(x) = \left( f(x),f'(x),\ldots,f^{(j)}(x), \ldots, f^{(n)}(x) \right) 
\] 
be the sequence of derivatives of $f(x)$.  
Let $V_f(x)$ be the number of sign changes  in the vector $D_f(x)$.  
The number of real zeros $\lambda$ of $f(x)$ 
(counting multiplicity) 
such that $a <\lambda \leq b$ is $V_f(a) - V_f(b) - 2\nu$ 
for some nonnegative integer $\nu$.  
\etNN

One proof of the theorem (Nathanson~\cite{nath22z}) 
uses the multiplicity vector $M_f(\lambda)$, 
and this leads to multiplicity matrices $M_f(\Lambda)$.  
For standard proofs of the Budan-Fourier theorem, see Basu, Pollack, and Roy~\cite{basu-poll-roy06}, Dickson~\cite{dick22}, and Jacobson~\cite{jaco74}.
An important early paper is Hurwitz~\cite{hurw12}.

\section{The multiplicity vector of a polynomial}

Let $f(x) \in K[x]$ be a polynomial of degree $n$ and let $f^{(j)}(x)$ 
denote the $j$th derivative of $f(x)$.   
For $j \in \{0,1,2,\ldots, n\}$, let $ \mu_j$ be the multiplicity of the element $\lambda \in K$
as a zero of the polynomial $f^{(j)}(x)$.  The vector 
\[
M_f(\lambda) = ( \mu_0,  \mu_1,\ldots,  \mu_n) 
\]  
is called the \emph{multiplicity vector of the polynomial $f(x)$ at $\lambda$}. 
The multiplicity vector $M_f(\lambda)$ is the zero vector if and only if $f^{(j)}(\lambda) \neq 0$ 
for all $j \in \{0,1,\ldots, n\}$.  

For example, the cubic polynomial $f(x) = x^3- 3x^2$ has sequence of derivatives 
\begin{align*}
f(x) & = x^3- 3x^2 = x^2(x-3) \\
f'(x) & = 3x^2 -6x= 3x(x-2)  \\
f''(x) & = 6x-6 = 6(x-1) \\
f^{(3)}(x) & =6 
\end{align*}
with multiplicity vectors 
\begin{align*}
M_f(0) & = (2,1,0,0) \\
M_f(1) & = (0, 0, 1,0) \\ 
M_f(2) & = (0,1,0,0)  \\
M_f(3) & = (1,0,0,0)  \\
M_f(\lambda) & = (0,0,0,0) \qquad \text{for all $\lambda \in K \setminus \{0,1,2,3\}$.}
\end{align*}

For the polynomial  $f(x) = 3x^5-7x^4+4x^3$ of degree 5 we have the sequence of derivatives 
\begin{align*}
f(x) & = 3x^5-7x^4+4x^3  = x^3  ( x-1)(3x-4 ) \\
f'(x) & =15x^4-28x^3 + 12 x^2  = x^2  \left(3x-2 \right)  \left(5x-6 \right) \\
f''(x) & = 60x^3 - 84 x^2 + 24 x  = 12x(x-1)(5x-2) \\
f^{(3)}(x) & = 180 x^2 - 168 x + 24   = 12 \left( 15x^2 - 14 x +2 \right) \\
f^{(4)}(x) & = 360 x - 168  = 24(15x-7) \\
f^{(5)}(x) & =   360.&
\end{align*}
The multiplicity vectors of $f(x)$ at $\lambda$ for  $\lambda \in \Lambda = ( 0,1)$ are 
\begin{align*}
M_f(0) & = (3,2,1,0,0,0) \\
M_f(1) & = (1,0,1,0,0, 0).
\end{align*}

A \emph{multiplicity vector}\index{multiplicity vector} 
is a vector 
\[
M = \left( \mu_0, \mu_1,\ldots, \mu_n \right)
\]
with nonnegative integer coordinates such that $\mu_n= 0$ 
and $\mu_j \geq 1$ implies $\mu_{j+1} = \mu_j -1$.

\bt                    \label{multiplicity:theorem:multiplicityVector} 
The multiplicity vector of a nonzero polynomial $f(x) \in K[x]$ at $\lambda \in K$ 
is a multiplicity vector.  
\et

\begin{proof}
If $f(x) = \sum_{j=0}^n c_jx^j$ has degree $n$, then $f^{(n)}(x) = n!c_n \neq 0$ 
and so $f^{(n)}(\lambda) \neq 0$ and $ \mu_n=0$. 
Equivalently, $\mu_j \geq 1$ implies $j \leq n-1$. 

If $\mu_j \geq 1$, then $j \leq n-1$ and there exists a polynomial $g_j(x)$ such that 
\[
f^{(j)}(x) = (x-\lambda)^{ \mu_j} g_j(x) \qqand g_j(\lambda) \neq 0.
\]
Differentiating this relation gives 
\begin{align*} 
f^{(j+1)}(x) & = { \mu_j}(x-\lambda)^{{ \mu_j}-1} g_j(x)  +  (x-\lambda)^{ \mu_j} g'_j(x) \\
& =  (x-\lambda)^{{ \mu_j}-1}   g_{j+1}(x)
\end{align*} 
where 
\[
 g_{j+1}(x) = { \mu_j} g_j(x)  +  (x-\lambda) g'_j(x) 
 \]
 and  
 \[
 g_{j+1}(\lambda) = { \mu_j} g_j(\lambda) \neq 0.
\]
Therefore, $ \mu_{j+1} =  \mu_j - 1$. 
This completes the proof. 
\end{proof}

\bl                                                 \label{multiplicity:lemma:Mf}
Let  $f(x) \in K[x]$  be a polynomial of degree $n$ and let $\lambda \in K$.    
The polynomial $f^{(j)}(x)$ 
has a zero of multiplicity $r$ at $x=\lambda$ if and only if 
$f^{(j+k)} (\lambda) = 0$ for all $k \in \{ 0,1,\ldots, r-1\}$ 
and $f^{(j + r)}(\lambda) \neq 0$. 
\el

\begin{proof} 
Let $\mu_{j+k}$ denote the multiplicity of $\lambda$ as a zero of $f^{(j+k)}(x)$ for 
all $k \in \{0,1,\ldots, n-j\}$.
By Theorem~\ref{multiplicity:theorem:multiplicityVector}, if $\mu_j = r$, 
then $\mu_{j+k} = r -k$ for $k \in \{0,1,\ldots, r\}$ 
and  so $f^{(j+k)} (\lambda) = 0$ for $k \in \{ 0,1,\ldots, r -1\}$ 
and $f^{(j + r)}(\lambda) \neq 0$. 

Conversely, if $f^{(j+k)} (\lambda) = 0$ for $k \in \{ 0,1,\ldots, r-1\}$ 
and $f^{(j + r)}(\lambda) \neq 0$, then 
$\mu_{j+k} \geq 1$ for $k \in \{ 0,1,\ldots, \mu_j -1\}$ 
and $\mu_{j+r} = 0$.  
Theorem~\ref{multiplicity:theorem:multiplicityVector} implies 
that if $\mu_{j+r} =0$ and 
$\mu_{j+r-1} \geq 1$, then $\mu_{j+r-1} = \mu_{j+r} +1 = 1$.  
Similarly, $\mu_{j+r-2} = \mu_{j+r-1} +1 = 2$.  Continuing inductively, we obtain 
$\mu_{j+r} = r$. This completes the proof. 
\end{proof}

\bl[Taylor's formula]                         \label{multiplicity:lemma:Taylor}
Let   $\lambda\in K$   and let $f(x) = \sum_{j=0}^n c_j (x-\lambda)^j \in K[x]$ 
be a polynomial of degree $n$.  
For all $j \in \{ 0,1,\ldots, n \}$, let $f^{(j)}(x)$ be the $j$th derivative of $f(x)$. Then 
\[
c_j = \frac{f^{(j)}(\lambda)}{j!}
\]
and $f^{(j)}(\lambda) = 0$ if and only if $c_j=0$.  
\el

\bt                                            \label{multiplicity:theorem:MultiplicityVector} 
Let $M = (\mu_0, \mu_1, \ldots, \mu_n)$ be a multiplicity vector. 
Let $\lambda \in K$ and let $f(x) = \sum_{j=0}^n c_j (x-\lambda)^j  \in K[x]$ 
be a polynomial of degree $n$.  
A necessary and sufficient condition that 
$M$ be the multiplicity vector 
of  $f(x)$ at $x=\lambda$ is that $c_j = 0$ if and only $\mu_j \geq 1$.
\et

\begin{proof}
Let $M$ be the multiplicity vector 
of  $f(x)$ at $x=\lambda$, that is, $M = M_f(\lambda)$.   
By Taylor's formula (Lemma~\ref{multiplicity:lemma:Taylor}), if  $\mu_j \geq 1$, 
then $f^{(j)}(\lambda) = 0$ and so $c_j = 0$. 
If  $\mu_j = 0$, then $f^{(j)}(\lambda) \neq 0$ and so $c_j \neq 0$.

Conversely, let $M = (\mu_0, \mu_1, \ldots, \mu_n)$ be a multiplicity vector 
and let $f(x) = \sum_{j=0}^n c_j (x-\lambda)^j$ be a polynomial of degree $n$ such that 
$c_j = 0$ if and only $\mu_j \geq 1$.  
If $\mu_j = 0$, then $f^{(j)}(\lambda) = j! c_j \neq 0$ and $f^{(j)}(x)$ 
has a zero of multiplicity $0 = \mu_j$ at $x=\lambda$. 
By Lemma~\ref{multiplicity:lemma:Mf}, if $\mu_j = r \geq 1$,  
then $\mu_{j+k}=  r-k \geq 1$ for $k \in \{ 0,1,\ldots, r-1 \}$ 
and $\mu_{j+r}=  0$.  
Therefore,  $c_{j+k} = 0$ for $k \in \{ 0,1,\ldots, r-1\}$ and $c_{j+r} \neq 0$.
Equivalently, by Taylor's formula, $f^{(j+k)}(\lambda)  = 0$ 
for $k \in \{ 0,1,\ldots, r-1\}$ and $f^{(j+r)}(\lambda) \neq 0$.
By Lemma~\ref{multiplicity:lemma:Mf}, the polynomial $f^{(j)}(x)$ 
has a zero of multiplicity $r = \mu_j$ at $x=\lambda$. 
This completes the proof. 
\end{proof}

\section{Multiplicity matrices}                 \label{multiplicity:section:MM}

Let $f(x)  \in K[x]$ be a polynomial of degree $n$ and let $f^{(j)}(x)$ be the $j$th 
derivative of $f(x)$.
Let 
\[
\Lambda = (\lambda_1, \lambda_2,\ldots, \lambda_m)
\]
 be a sequence of distinct elements of $K$.  
For $i \in \{1, 2, \ldots, m \}$, 
let $ \mu_{i,j}$ be the multiplicity of $\lambda_i$ as a zero of the polynomial $f^{(j)}(x)$. 
The \emph{multiplicity matrix of $f(x)$ with respect to $\Lambda$}\index{multiplicity matrix} 
is the $m \times (n+1)$ matrix 
\beq        \label{multiplicity:MultiplicityMatrix}
M_f(\Lambda) =  
\bmat 
 \mu_{1,0} &  \mu_{1,1} &  \ldots & \mu_{1,j} & \ldots &  \mu_{1,n} \\
 \mu_{2,0} &  \mu_{2,1} &  \ldots & \mu_{2,j} & \ldots &  \mu_{2,n} \\
\vdots  &  & & & \vdots  \\
 \mu_{i,0} &  \mu_{i,1} &  \ldots & \mu_{i,j} & \ldots &  \mu_{i,n} \\
\vdots  &  & & & \vdots  \\
 \mu_{m,0} &  \mu_{m,1} &  \ldots & \mu_{m,j} &\ldots &  \mu_{m,n}   
\emat. 
\eeq
For all $i \in \{1,2,\ldots,m\}$, the $i$th row of $M_f(\Lambda)$ is the 
multiplicity vector $M_f(\lambda_i)$. 

For example, the multiplicity matrix of the polynomial $f(x) = x^4 -4x^3 + 4x$ and 
the sequence $\Lambda = (0,1, 2 )$ is 
\[
M_f(\Lambda) =  
\bmat 
1 & 0 & 1 & 0 & 0 \\
0 & 0 &  0 & 1 & 0 \\
0 & 0 & 1 & 0 & 0 
\emat.
\]

The $j$th column sum of the $m \times (n+1)$ matrix $M = \bmat m_{i,j} \emat$ is
\[
\colsum_j(M) = \sum_{i=1}^{m} \mu_{i,j}.
\]
The definition of the multiplicity  integer $\mu_{i,j}$ implies that for all $j \in \{1,2,\ldots, m \}$  
there exists a nonzero polynomial $g_j(x)$ such that 
\[
f^{(j)}(x) = g_j(x) \prod_{i=1}^{m} (x-\lambda_i)^{\mu_{i,j}} 
\]
and
\[
g_j(\lambda_i) \neq 0 \quad  \text{ for all } \lambda_i \in \Lambda.
\]
Because $\deg\left( f^{(j)}(x) \right)  = \deg(f(x)) - j = n - j $, we have 
the following constraint on the column sums of the multiplicity matrix: 
\begin{align*}
 \colsum_j\left( M_f(\Lambda) \right)  
 = \sum_{i=1}^{m} \mu_{i,j} & \leq   \sum_{i=1}^{m} \mu_{i,j} +  \deg\left( g_j(x) \right) \\
&  =  \deg\left( f^{(j)}(x)  \right)  = n - j 
\end{align*}
for all $j \in \{0,1,\ldots, n\}$. 

For every nonzero polynomial $f(x)$ and every nonzero element $c$, 
we have $M_f(\Lambda) = M_{cf}(\Lambda)$.  
Thus, the multiplicity matrix of a polynomial is always the multiplicity matrix of a monic polynomial.

An \emph{$m\times (n+1)$ multiplicity matrix}\index{multiplicity matrix} 
is an $m \times (n+1)$ matrix  
\[ 
M = \bmat 
 \mu_{1,0} &  \mu_{1,1} &  \mu_{1,2} & \ldots &  \mu_{1,n} \\
 \mu_{2,0} &  \mu_{2,1} &  \mu_{2,2} & \ldots &  \mu_{2,n} \\
\vdots  &  & & & \vdots  \\
 \mu_{m,0} &  \mu_{m,1} &  \mu_{m,2} &\ldots &  \mu_{m,n}   
\emat
\]
with rows enumerated by $i \in \{ 1,2,\ldots, m \}$ 
and columns enumerated by $j \in \{ 0,1,\ldots, n \}$, 
such that 
\benum
\item[(i)]
every row of $M$ is a multiplicity vector, and   
\item[(ii)] 
 for all $j \in \{0,1,\ldots, n\}$, 
 \[
 \colsum_j(M) =  \sum_{i=1}^{m} \mu_{i,j} \leq   n - j.
 \] 
\eenum
Every matrix obtained from a multiplicity matrix by a permutation of the rows 
is also a multiplicity matrix.

For example, $\bmat 2 &1 & 0 & 0 \\ 0 & 1 & 0 & 0\emat$ is a $(2,4)$ multiplicity matrix but 
 $\bmat 2 &1 & 0  \\ 0 & 1 & 0 \emat$ is not a $(2,3)$ multiplicity matrix.

The study of multiplicity matrices is the study of interactions between the zeros of a polynomial.  
We ask the following questions (cf. Open Problems~\eqref{Multiplicity:problem:1}, 
~\eqref{Multiplicity:problem:2},  
and~\eqref{Multiplicity:problem:3} in Section~\ref{multiplicity:section:OpenProblems}) . 

Given an $m \times (n+1)$ multiplicity matrix $M$ and a sequence 
 $\Lambda = (\lambda_1, \lambda_2, \ldots, \lambda_m)$ 
 of elements of $K$, does there exist a polynomial $f(x)$ of degree $n$ 
such that $M_f(\Lambda) = M$? 

Given an $m \times (n+1)$ multiplicity matrix $M$, 
do there exist a sequence 
 $\Lambda = (\lambda_1, \lambda_2, \ldots, \lambda_m)$  
 of elements of $K$ and  a polynomial $f(x)$ of degree $n$ 
such that $M_f(\Lambda) = M$? 

Describe the multiplicity matrices that are not multiplicity matrices of polynomials.

\bex
The multiplicity matrix 
\[
M = \bmat 2 &1 & 0 & 0 \\ 0 & 1 & 0 & 0\emat 
\]
is the multiplicity matrix of the cubic polynomial 
\[
f(x) = x^3 -\frac{3}{2}x^2 = (x-1)^3 + \frac{3}{2} (x-1)^2 - \frac{1}{2} 
\]
with respect to the sequence  $\Lambda = ( 0,1)$. 
Moreover, $f(x)$ is the unique monic cubic polynomial such that $M = M_f(\Lambda)$. 
\eex

\begin{proof} 
If   $M$ is the multiplicity matrix of a polynomial, then the polynomial is cubic.  
The  matrix $M$ is the multiplicity matrix of a monic cubic polynomial $f(x)$ 
with respect to  $\Lambda$ 
if and only if the multiplicity vector of $f(x)$ at $x = 0$ is 
\beq                        \label{multiplicity:MultiplicityExample1-0}
M_f(0) = \bmat 2 &1 & 0 & 0  \emat 
\eeq
and the multiplicity vector of $f(x)$ at $x = 1$ is 
\beq                        \label{multiplicity:MultiplicityExample1-1}
M_f(1) = \bmat 0 & 1 & 0 & 0  \emat. 
\eeq

By Theorem~\ref{multiplicity:theorem:MultiplicityVector}, 
condition~\eqref{multiplicity:MultiplicityExample1-0} implies that there is a  nonzero element 
$a \in K$ such that 
\[
f(x) = x^3 + ax^2.
\] 
Condition~\eqref{multiplicity:MultiplicityExample1-1} implies that there are nonzero elements 
$b$ and $c $  such that 
\begin{align*}
f(x) & = (x-1)^3 + b(x-1)^2 + c  \\
& = x^3 +(b-3)x^2 + (3-2b)x+(b+c-1) .
\end{align*}
Solving the equation 
\[
x^3 + ax^2 = x^3 +(b-3)x^2 + (3-2b)x+(b+c-1) 
\]
we obtain  
\begin{align*}
a & = b-3 \\
0 & = 3-2b \\
0 & = b+c-1
\end{align*}
and so 
\[
a = -\frac{3}{2}, \qquad b= \frac{3}{2}, \qquad c= -\frac{1}{2}.  
\]
Thus, 
\[
f(x) = x^3 -\frac{3}{2}x^2 = (x-1)^3 + \frac{3}{2} (x-1)^2 - \frac{1}{2}
\]
is the unique monic cubic polynomial such that $M = M_f(\Lambda)$. 
\end{proof}

\bex  
The  multiplicity matrix 
\[
M' = \bmat 3 & 2 &1 & 0 & 0 \\ 0 & 1 & 0 & 1& 0\emat 
\]    
is not the multiplicity matrix of a polynomial 
with respect to the sequence  $\Lambda = ( 0,1)$.  
\eex

\begin{proof} 

The multiplicity matrix $M'$ is the multiplicity matrix of a monic quartic   
polynomial $f(x)$ 
with respect to the sequence $\Lambda$ 
if and only if the multiplicity vector of $f(x)$ at $x = 0$ is 
\beq                        \label{multiplicity:MultiplicityExample2-0}
M'_f(0) = \bmat 3 & 2 &1 & 0 & 0  \emat 
\eeq
and the multiplicity vector of $f(x)$ at $x = 1$ is 
\beq                        \label{multiplicity:MultiplicityExample2-1}
M'_f(1) = \bmat 0 & 1 & 0 & 1 & 0  \emat. 
\eeq
By Theorem~\ref{multiplicity:theorem:MultiplicityVector}, 
condition~\eqref{multiplicity:MultiplicityExample2-0} implies that there is a nonzero element  
$a$  such that 
\[
f(x) = x^4 + ax^3.
\] 
Condition~\eqref{multiplicity:MultiplicityExample2-1} implies that there are nonzero elements 
$b$ and $c $  such that 
\begin{align*}
f(x) & = (x-1)^4 + b(x-1)^2 + c  \\
& = x^4 -4x^3+(b+6)x^2 -(2b+4)x + (b+c+1).
\end{align*} 
Therefore, 
\[
x^4 + ax^3 = x^4 -4x^3+(b+6)x^2 -(2b+4)x + (b+c+1).
\]
Comparing coefficients of the quadratic and  linear powers of $x$, we obtain 
\[
b+6 = 0 \qqand 2b+4 = 0
\]
and so $b = -6 = -2$, which is absurd.  Therefore, the multiplicity matrix $M'$ 
is not the multiplicity matrix of a polynomial.  
\end{proof}

\section{Solving the equation $M = M_f(\Lambda)$}

Given an $m \times (n+1)$ multiplicity matrix $M$, we would like to find 
all monic polynomials $f(x) \in K[x]$ of degree $n$ 
and sequences $\Lambda$ of length $m$ of elements of $K$ such that 
$M = M_f(\Lambda)$.  The following results describe ``equivalence classes'' 
of solutions of $M = M_f(\Lambda)$.

Let $j$ and $k$ be  positive integers.   
The \emph{$j$th falling factorial}\index{falling factorial} is the polynomial 
\[
(x)_j = x(x-1)\cdots (x-j+1). 
\]
Let $(x)_0 = 1$.  We have $(k)_j= 0$ if $j > k$.
If $h(x) = (x-\lambda)^k$, then $h^{(j)}(x) = (k)_j (x-\lambda)^{k-j}$.

An \emph{automorphism}\index{automorphism} of the field $K$ 
is a function $\sigma:K \rightarrow K$ such that 
\benum
\item
$\sigma(a+b) = \sigma(a) + \sigma(b)$ for all $a,b \in K$,
\item
$\sigma(ab) = \sigma(a)  \sigma(b)$ for all $a,b \in K$,
\item
$\sigma(1) = 1$.
\eenum
If $\sigma$ is an automorphism of $K$, then $\sigma(k) = k$ for all integers $k$ 
and so $\sigma \left( (k)_j \right) =  (k)_j $ for all nonnegative integers $j$ and $k$.

 \bt                   \label{multiplicity:theorem:sigma}
Let $\sigma$ be an automorphism of the field $K$.  
Let
\[
\Lambda = \left( \lambda_1, \lambda_2, \ldots, \lambda_m \right)
\]
be a sequence of distinct elements of $K$ and let 
\[
\Lambda_{\sigma} = ( \kappa_1,  \kappa_2, \ldots,  \kappa_{m} ) 
\]
where 
\[
\kappa_i = \sigma(\lambda_i) 
\]  
for all $i \in \{1,2,\ldots, m\}$.  
Let 
\[
f(x) = \sum_{k=0}^n c_k x^k \in K[x]
\]
be a polynomial of degree $n$ and let  
\[
f_{\sigma}(x) = \sum_{k=0}^n \sigma( c_k) x^k \in K[x]. 
\]  
Let $M = \bmat \mu_{i,j} \emat$ be an $m \times (n+1)$ multiplicity matrix.
Then  $M = M_f(\Lambda)$ if and only if $M = M_{f_{\sigma}} (\Lambda_{\sigma})$. 
\et

\begin{proof} 
Because $\left( \Lambda_{\sigma} \right)_{\sigma^{-1}} = \Lambda$ and 
$\left(f_{\sigma}\right)_{\sigma^{-1}} (x) = f(x)$, 
it suffices to prove that $M = M_f(\Lambda)$ implies 
$M = M_{f_{\sigma}} (\Lambda_{\sigma})$. 

The $j$th derivative of $f_{\sigma}(x)$ is  
\[
f_{\sigma}^{(j)}(x) =  \sum_{k=j}^n (k)_j \sigma( c_k) x^{k-j}. 
\] 
For all $\kappa_i \in \Lambda_{\sigma}$,  we have 
\begin{align*}
f_{\sigma}^{(j)}(\kappa_i) & = f_{\sigma}^{(j)}( \sigma(\lambda_i)) =  \sum_{k=j}^n (k)_j \sigma( c_k) \sigma(\lambda_i)^{k-j} \\
& = \sigma\left( \sum_{k=j}^n (k)_j   c_k \lambda_i^{k-j}  \right)  = \sigma\left( f^{(j)}  (\lambda_i)   \right).
\end{align*}
Because $\sigma$ is an automorphism,
for all $i \in \{1,2,\ldots, m\}$ and $j \in \{0,1,\ldots, n\}$ we have 
$f_{\sigma}^{(j)}(\kappa_i) =0$ if and only if $ f^{(j)}  (\lambda_i) = 0$.  
By Lemma~\ref{multiplicity:lemma:Mf}, the multiplicity of $\kappa_i$ as a zero of  
$f_{\sigma}^{(j)}(x)$ equals  the multiplicity of $\lambda_i$ as a zero of  
$f^{(j)}(x)$, and so  
$M_{f_{\sigma}}(\Lambda_{\sigma}) = M_{f}(\Lambda)$.   
This completes the proof.  
\end{proof}

Sequences $\Lambda =\left( \lambda_1, \lambda_2, \ldots, \lambda_m \right)$ 
and $\Lambda' = ( \kappa_1,  \kappa_2, \ldots,  \kappa_{m} )$ in $K$ 
are \emph{affine equivalent}\index{affine equivalent} 
if there exist $r, s \in K$ with $r \neq 0$ such that 
\[
\lambda_i = r  \kappa_i  + s
\] 
for all $i \in \{1,2,\ldots, m\}$.

 \bt                   \label{multiplicity:theorem:affine}
Let $r,s \in K$ with $r \neq 0$.  
Let
\[
\Lambda = \left( \lambda_1, \lambda_2, \ldots, \lambda_m \right)
\]
be a sequence of distinct elements of $K$ and let 
\[
\Lambda_{r,s} = ( \kappa_1,  \kappa_2, \ldots,  \kappa_{m} ) 
\]
where 
\[
\lambda_i = r \kappa_i + s 
\]  
for all $i \in \{1,2,\ldots, m\}$.  
Let $f(x) \in K[x]$ 
be a polynomial of degree $n$ and let  
\[
f_{r,s}(x) = f(rx+s)  \in K[x]. 
\]  
Let $M = \bmat \mu_{i,j} \emat$ be an $m \times (n+1)$ multiplicity matrix.
Then  $M = M_f(\Lambda)$ if and only if $M = M_{f_{r,s}} \left( \Lambda_{r,s} \right)$. 
\et

\begin{proof} 
Let $f(x)$ be a polynomial of degree $n$ such that $M = M_{f}(\Lambda)$. 
For all $j \in \{0,1,\ldots, n\}$ there is a polynomial $g(x)$ of degree 
$n-j - \sum_{i=1}^{m} \mu_{i,j}$ such that 
\[
f^{(j)}(x) = g_{j}(x) \prod_{i=1}^{m}   (x- \lambda_i)^{\mu_{i,j}} 
\]
and 
\[
g_{j}(\lambda_i)  \neq 0
\] 
for all $i \in \{1,2,\ldots, m\}$.

The polynomial $f_{r,s}(x) = f\left( rx+s \right)$ has degree $n$ and 
\[
f_{r,s}(\kappa_i) = f\left( r\kappa_i +s \right)  = f(\lambda_i) 
\]
for all $i \in \{ 1,2, \ldots, m \}$.
Moreover, for all $j \in \{ 0,1, \ldots, n\}$ we have 
\begin{align*}
f_{r,s}^{(j)}(x) 
& = r^j f^{(j)}\left( rx+s \right) \\ 
& = r^j g_j \left( rx+s  \right) 
\prod_{i=1}^{m} \left( rx+s - \lambda_i \right)^{\mu_{i,j}} \\
& = r^j g_j \left(rx+s \right) 
\prod_{i=1}^{m} \left( r(x-\kappa_i) \right)^{\mu_{i,j}} \\
& = r^{j + \sum_{i=1}^{m} \mu_{i,j}}  g_j \left(rx+s \right) 
\prod_{i=1}^{m} \left( x-\kappa_i \right)^{\mu_{i,j}} \\
 & =   g_{r,s,j}\left( x  \right) 
  \prod_{i=1}^{m}   \left( x-\kappa_i \right)^{\mu_{i,j}} \\
\end{align*} 
where 
\[
g_{r,s,j} (x) = r^{j + \sum_{i=1}^{m} \mu_{i,j}} 
 g_j \left(  rx+s \right) 
 \]
and 
\[
g_{r,s,j}( \kappa_i) =  r^{j + \sum_{i=1}^{m} \mu_{i,j}}  
 g_j \left( r \kappa_i +s\right)
 =  r^{j + \sum_{i=1}^{m} \mu_{i,j}}  
 g_j \left( \lambda_i \right)  \neq 0 
\]
for all $i \in \{ 1,2, \ldots, m\}$. 
Therefore, $M = M_{f_{r,s}}(\Lambda_{r,s})$.  

Let $r' = 1/r$ and $s' = -s/r$.  The observation that 
$\lambda_i = r \kappa_i + s$ if and only if $\kappa_i =  r' \lambda_i + s'$ 
and so 
\[
\left(f_{r,s}\right)_{r',s'} (x) = f(x) 
\]
and 
\[
\left( \Lambda_{r,s}\right)_{r',s'} = \Lambda 
\]
 completes the proof. 
\end{proof}

\bt
Let $\sigma$ be an automorphism of $K$ 
and let $r,s \in K$ with $r \neq 0$.  
Let
\[
\Lambda = \left( \lambda_1, \lambda_2, \ldots, \lambda_m \right)
\]
be a sequence of elements of $K$ and let 
\[
\tilde{\Lambda} = \left( \tilde{\lambda}_1, \tilde{ \lambda}_2, \ldots, \tilde{ \lambda}_{m}  \right)
\]
where
\[
\tilde{\lambda_i }= r\sigma(\lambda_i)+s
\]
for all $i \in \{1,2,\ldots, m\}$. 
Let $M = \bmat \mu_{i,j} \emat$ be an $m \times (n+1)$ multiplicity matrix. 
There is a polynomial $f(x) \in K[x]$ such that $M = M_f(\Lambda)$ 
if and only if there is a  polynomial $\tilde{f}(x) \in K[x]$ such that 
$M = M_{ \tilde{f} } \left( \tilde{\Lambda} \right)$. 
\et

\begin{proof}
This follows from Theorems~\ref{multiplicity:theorem:sigma} 
and~\ref{multiplicity:theorem:affine} and the observation that 
$\tilde{\Lambda} = \left(\Lambda_{\sigma}\right)_{r,s}$.    
\end{proof}

\bt                              \label{multiplicity:theorem:affine-01} 
Let $m \geq 1$ and let $M$ be an $m \times (n+1)$ multiplicity matrix. 
There exists a sequence $\Lambda'$ of distinct elements of $K$ and a polynomial $g(x) \in K[x]$ 
such that $M = M_g(\Lambda')$ if and only if there exists a sequence $\Lambda$ 
of distinct elements of $K$ with $0,1 \in \Lambda$ and a polynomial $f(x) \in K[x]$ 
such that $M = M_{f}(\Lambda)$.
\et

\begin{proof}
Let $\Lambda' = (\kappa_1, \kappa_2,\ldots, \kappa_{m})$ be a 
sequence of distinct elements of $K$.  For all $i \in \{1,2,\ldots, m\}$, let
\[
\lambda_i = \frac{\kappa_i - \kappa_1}{\kappa_2 - \kappa_1}
\]
and let 
\[
\Lambda = (\lambda_1, \lambda_2,\ldots, \lambda_m).
\]
The sequences $\Lambda'$ and $\Lambda$ are affine equivalent 
and $\lambda_1 = 0$ and $\lambda_2 = 1$. 
This completes the proof. 
\end{proof} 

\bc                              \label{multiplicity:corollary:affine} 
Let $M$ be a $2 \times (n+1)$ multiplicity matrix. 
Let 
\[
\Lambda_0 = ( 0,1 )
\]
and let
\[
\Lambda = ( \lambda_1, \lambda_2 ) \subseteq K
\]
with $\lambda_1 \neq \lambda_2$.  
There exists a polynomial $f_0(x)$  such that $M = M_{f_0}(\Lambda_0)$  
if and only if there exists a polynomial $f(x)$  such that $M = M_{f}(\Lambda)$. 
\ec

Thus, in the study of multiplicity matrices $M$ with two rows, it suffices to consider 
only polynomials $f(x)$ such that $M = M_f(\Lambda)$ 
with respect to the sequence $\Lambda_0 = (0,1)$. 
The story is more complicated  for multiplicity matrices $M$ with three or more rows, 
that is, for sequences $\Lambda$ that contain three or more elements.  
In an earlier version of this paper, I had asked the following question:  
Let $\Lambda_1 = ( 0,1,2 )$ and  $\Lambda_2 =  ( 0,1,3 )$.   
Construct a  $3 \times (n+1)$  multiplicity matrix $M$ with the following property: 
There exists a polynomial $f(x)$ of degree $n$ such that $M = M_f(\Lambda_1)$ 
but there does not 
exist a polynomial $g(x)$ of degree $n$ such that $M = M_g(\Lambda_2)$.   

Independently, Sergei Konyagin and Noah Kravitz (personal communications)  
observed that if  the $3 \times 3$ multiplicity matrix 
\[
M = \bmat 1 & 0 & 0 \\ 0 & 1 & 0 \\ 1 & 0 & 0 \emat
\]
were the multiplicity matrix of a monic quadratic polynomial $f(x)$ for the sequence 
$\Lambda = ( 0, 1,\lambda )$ 
for some $\lambda \notin \{ 0,1 \}$, then 
\[
f(x) = x^2 - 2x = x(x-\lambda)  
\]
and so $\lambda = 2$.

\section{Multiplicity matrices with no polynomial}
For every nonnegative integer $j$ and $m$-tuple   of nonnegative integers
$(k_1, k_2, \ldots, k_{m})$ such that $\sum_{i=1}^{m} k_i = j$, 
we define the multinomial coefficient 
\[
\binom{j}{k_1,k_2, \ldots, k_{m}} = \frac{j!}{k_1! k_2! \cdots k_{m}!}.  
\]

We recall  the \emph{Leibniz rule}\index{Leibniz rule} 
for derivatives of products of finitely many functions.  

\bl[Leibniz]         \label{multiplicity:lemma:Leibniz} 
Let $h_i(x) \in K[x]$ for all $i \in \{1,2,\ldots, m \}$ and let $f(x) = \prod_{i=1}^{m} h_i(x)$.  
For all nonnegative integers $j$, 
\[
f^{(j)}(x) = \sum_{(k_1, k_2, \ldots, k_{m})} \binom{j}{k_1, k_2, \ldots, k_{m}} \prod_{i=1}^{m} h_i^{(k_i)}(x) 
\]
where the summation is over all $m$-tuples of nonnegative integers 
such that $\sum_{i=1}^{m} k_i = j$ 
and $ \binom{j}{k_1, k_2,\ldots, k_{m}}$ is the multinomial coefficient.  
\el

\begin{proof}
The proof is by induction on $m$ for $m \geq 2$.   The starting case $m=2$ is the 
usual Leibniz rule for the derivatives of a product of two functions.  
\end{proof}

To understand interactions of zeros of polynomials, it is useful  
to construct multiplicity matrices that are not multiplicity matrices of polynomials.    
We begin with a simple observation.  
Let $M = \bmat \mu_{i,j} \emat$ be an $m\times (n+1)$  matrix.   
For $\ell \in \{0,1,\ldots, n \}$,  the 
\emph{$\ell$-truncation}\index{truncation} of $M$, 
denoted $M^{(\ell)}$, is the $m \times (n+1-\ell)$ submatrix 
of $M$ obtained by deleting the first $\ell$ columns of $M$.  
Note that $M^{(0)} = M$.  
For example, the 2-truncation of 
\[
M = \bmat 5 & 4 & 3 & 2 & 1 & 0 \\ 0 & 2 & 1 & 0 & 0 & 0 \emat
\]
is
\[
M^{(2)} = \bmat 3 & 2 & 1 & 0 \\  1 & 0 & 0 & 0 \emat.
\]

\bl                               \label{multiplicity:lemma:truncate}
Let $M = \bmat \mu_{i,j} \emat$ be an $m\times (n+1)$ multiplicity matrix 
and let $f(x) \in K[x]$ be a polynomial of degree $n$.  
For $\ell \in \{0,1,\ldots, n \}$, let $M^{(\ell)}$ be the $\ell$-truncation of $M$ 
and let $f^{(\ell)}(x)$ be the $\ell$th derivative of $f(x)$.
Let $\Lambda$ be a sequence of elements of $K$.  
If $M = M_f(\Lambda)$, then $M^{(\ell)} =M_{f^{(\ell)}}(\Lambda)$.  
\el

\begin{proof}  
For all $i \in \{1,2,\ldots , m \}$ and $j \in \{0,1,\ldots, n- \ell \}$,  
if  $\lambda_i$ is a zero of $f^{(\ell)}(x)$ of multiplicity $\mu^{(\ell)}_{i,j}$, 
then 
\[
f^{(\ell+j)}(x) = \left( f^{(\ell)} \right)^{(j)}(x) = g(x)(x-\lambda_i)^{ \mu^{(\ell)}_{i,j} } 
\]   
for some polynomial $g(x)$ with $g(\lambda_i) \neq 0$, and so $\mu^{(\ell)}_{i,j} = \mu_{\ell+j}$. 
This completes the proof. 
\end{proof}

\bt                              \label{multiplicity:theorem:truncate}
Let $M' = \bmat \mu_{i,j} \emat$ be an $m\times (n'+1)$ multiplicity matrix 
that is not the multiplicity matrix of a polynomial. 
Let $n = n' + \ell$.  If $M = \bmat \mu_{i,j} \emat$ 
is an $m\times (n+1)$ multiplicity matrix 
such that $M'$ is the $\ell$-truncation of $M$, then $M$ is  
not the multiplicity matrix of a polynomial. .
\et

\begin{proof}
This follows from Lemma~\ref{multiplicity:lemma:truncate}.  
\end{proof}

\bt                         \label{multiplicity:theorem:obstruction-1} 
The $2 \times 4$ multiplicity matrix 
\[ 
M = \bmat 
 2 & 1 & 0 & 0 \\
 1 & 0  &1 & 0   
\emat
\]
is not the multiplicity matrix of a polynomial.  
\et

\begin{proof}   
By Corollary~\ref{multiplicity:corollary:affine}, the matrix $M$ is the multiplicity matrix 
of a polynomial with respect to some sequence $\Lambda = (\lambda_1, \lambda_2)$ 
if and only if $M$ is the multiplicity matrix of a polynomial 
with respect to the sequence $\Lambda_0 = (0,1)$. 

By Theorem~\ref{multiplicity:theorem:MultiplicityVector}, 
if $f(x)$ is a monic cubic polynomial such that $M = M_f(\Lambda)$, 
then the first and second rows of $M$ imply that there are nonzero 
elements $a, b \in K$ such that 
\begin{align*}
f(x) & =    x^3 + ax^2  =  (x-1)^3 + b(x-1)  \\
& =  x^3 - 3x^2 + (b+3)x - (b +1)  
\end{align*}
and so 
\[
b= - 3 = - 1
\]
which is absurd.  This completes the proof. 
\end{proof}

\bt                         \label{multiplicity:theorem:obstruction-2} 
Let  $M$ be a $2 \times (n+1)$ multiplicity matrix of the form 
\[ 
M = \bmat 
 \mu_{1,0} &  \mu_{1,1} &  \mu_{1,2} & \ldots &  \mu_{1,n-4} & 2 & 1 & 0 & 0 \\
 \mu_{2,0} &  \mu_{2,1} &  \mu_{2,2} & \ldots &  \mu_{2,n-4} & 1 & 0  &1 & 0   
\emat.
\]
There exists no polynomial $f(x)$ and no sequence $\Lambda = (\lambda_1, \lambda_2)$ 
such that $M = M_f(\Lambda)$. 
\et

\begin{proof}
This follows from Theorems~\ref{multiplicity:theorem:truncate} 
and~\ref{multiplicity:theorem:obstruction-1}. 
\end{proof}

\bl                      \label{multiplicity:lemma:ApplyLeibniz} 
Let $p$ be a positive integer, let $\lambda \in K$, 
let $g(x) \in K[x]$ be a polynomial such that $g(\lambda) \neq 0$, and let
 \[
 f(x) = (x - \lambda)^p g(x).
 \]
Then
\begin{align*}
f^{(j)}(\lambda) & = 0  \quad\qquad\qquad\qquad \text{ for $0 \leq j \leq p-1 $} \\
f^{(p)}(\lambda) & \neq 0  \\
f^{(j)}(\lambda) &= (j)_p g^{(j-p)} (\lambda) \qquad \text{ for $j \geq p+1$} . 
\end{align*}

\el

\begin{proof}
For $j \in \{0, 1, \ldots, p-1\}$, the Leibniz rule (Lemma~\ref{multiplicity:lemma:Leibniz})  
gives  
\begin{align*}
f^{(j)}(x) & = \sum_{k=0}^j \binom{j}{k} (p)_k (x - \lambda)^{p-k} g^{(j-k)}(x) \\
& =  (x - \lambda) \sum_{k=0}^j \binom{j}{k} (p)_k (x - \lambda)^{p-1-k} g^{(j-k)}(x) 
\end{align*}
and so $f^{(j)}(\lambda) = 0$.

If $j  = p$, then 
\begin{align*}
f^{(p)}(x) 
& = \sum_{k=0}^p \binom{p}{k} (p)_k (x - \lambda)^{p-k} g^{(p-k)}(x) \\ 
& =  (x - \lambda) \sum_{k=0}^{p-1} \binom{p}{k} (p)_k (x - \lambda)^{p-1-k} g(x) 
+ p! g(x) 
\end{align*} 
and $f^{(p)}(\lambda) = p! g(\lambda) \neq 0$. 

If $j  \geq p+1$, then $(p)_k =0$ for $k > p$ implies 
\begin{align*}
f^{(j)}(x) 
& = \sum_{k=0}^j \binom{j}{k} (p)_k (x - \lambda)^{p-k} g^{(j-k)}(x) \\ 
& =  \sum_{k=0}^p \binom{j}{k} (p)_k (x - \lambda)^{p-k} g^{(j-k)}(x) \\ 
& =  (x - \lambda) \sum_{k=0}^{p-1} \binom{j}{k} (p)_k (x - \lambda)^{p-1-k} g^{(j-k)}(x) 
+\binom{j}{p}  p! g^{(j-p)}(x) 
\end{align*} 
and so 
\[
f^{(j)}(\lambda) = \binom{j}{p}  p! g^{(j-p)} (\lambda) =  (j)_p g^{(j-p)} (\lambda).
\] 
This completes the proof. 
\end{proof}

\bt                         \label{multiplicity:theorem:obstruction-3} 
Let $M$ be a $2 \times (n+1)$ multiplicity matrix 
\[ 
M = \bmat 
 \mu_{1,0} &  \mu_{1,1} &  \mu_{1,2} & \ldots &  \mu_{1,n-1} &  \mu_{1,n} \\
 \mu_{2,0} &  \mu_{2,1} &  \mu_{2,2} & \ldots &  \mu_{2,n-1}  &  \mu_{2,n}   
\emat
\]
such that 
\[
\emph{colsum}_0(M) =  \mu_{1,0}  + \mu_{2,0} = n.
\]
Let
\[
p =  \mu_{1,0}  \qqand q =  \mu_{2,0}.
\]
The matrix $M$ is the multiplicity matrix of a monic polynomial $f(x)$ of degree $n$ 
with respect to $\Lambda_0 = (0,1)$ if and only if 
\beq                      \label{multiplicity:obstruct-f} 
f(x) =  x^{p}  (x- 1)^{q}
\eeq
and 
\beq                      \label{multiplicity:obstruct-M} 
M = \bmat 
p &  p-1 & p-2& \ldots & 2 & 1 & 0 & 0 & \cdots & 0 \\
q &  q-1 & q-2&  \ldots & 1 & 0 & 0 & 0 & \cdots   & 0
\emat,
\eeq
that is, 
\beq                      \label{multiplicity:obstruct-M-p} 
\mu_{1,j} = 
\begin{cases}
p - j & \text{ for  $j \in \{0,1,\ldots, p - 1 \}$ } \\
0 & \text{for   $j \in \{ p, p+1,  \ldots, n\}$} 
\end{cases}
\eeq 
and
\beq                      \label{multiplicity:obstruct-M-q} 
\mu_{2,j} = 
\begin{cases}
q - j & \text{ for $j \in \{0,1,\ldots, q- 1 \}$ } \\
0 & \text{for   $j \in \{q, q+1,  \ldots, n\}$}.
\end{cases}
\eeq 
\et

For example, if $p = 3$, $q = 2$, and $n = 5$, 
then $f(x) = x^3 (x-1)^2$ and 
\[
M = M_f(\Lambda_0) = \bmat 
3 & 2  & 1 & 0 & 0 & 0 \\
2  & 1 & 0 & 0 & 0 & 0 
\emat. 
\]

\begin{proof}
By Lemma~\ref{multiplicity:lemma:ApplyLeibniz}, 
the multiplicity matrix of the polynomial $f(x) =  x^{p}  (x- 1)^{q}$ 
with respect to $\Lambda_0 = (0,1)$ 
is the $2 \times (n+1)$ matrix $M$ defined by~\eqref{multiplicity:obstruct-M}.  

Conversely, let $f(x)$ be a monic polynomial of degree $n$ such that $M = M_f(\Lambda_0)$.  
Because $f(x)$ has a zero at $x = 0$ of multiplicity $p$ 
and a zero at $x = 1$ of multiplicity $q$  
and because $p  + q = n$ 
it follows that 
\[
f(x) =  x^{p}   (x- 1)^{q}.
\]
The multiplicity relations~\eqref{multiplicity:obstruct-M-p} 
and~\eqref{multiplicity:obstruct-M-q} follow from Lemma~\ref{multiplicity:lemma:ApplyLeibniz}. 
This completes the proof.  
\end{proof}

Theorem~\ref{multiplicity:theorem:obstruction-3} enables us to construct 
multiplicity matrices that are not multiplicity matrices of polynomials. 
For example, there are six equivalence classes of 
$2 \times 6$ multiplicity matrices  such that 
$\mu_{0,0} = 3$ and $\mu_{1,0} = 2$: 
\[
M_1 = \bmat
3 & 2 & 1 & 0 & 0 & 0 \\
2 & 1 & 0 & 0 & 0 & 0
\emat,
\qquad 
M_2 = \bmat
3 & 2 & 1 & 0 & 0 & 0 \\
2 & 1 & 0 & 0 & 1 & 0
\emat,
\]
\[
M_3 = \bmat
3 & 2 & 1 &  0 & 0  & 0\\
2 & 1 & 0 & 1 & 0 & 0
\emat, 
\qquad 
M_4 = \bmat
3 & 2 & 1 &  0 & 0  & 0\\
2 & 1 & 0 & 2 & 1 & 0
\emat,
\]
\[
M_5 = \bmat
3 & 2 & 1 & 0 & 1 & 0 \\
2 & 1 & 0 & 0 & 0 & 0
\emat, 
\qquad 
M_6 = \bmat
3 & 2 & 1 & 0 & 1 & 0 \\
2 & 1 & 0 & 1 & 0 & 0
\emat. 
\]
By Theorem~\ref{multiplicity:theorem:obstruction-3}, the matrix $M_1$ 
is the unique matrix in this list 
that is the multiplicity matrix of a polynomial, 
and, with $\Lambda_0 = ( 0,1)$,  the unique monic quintic polynomial $f(x)$ such that  $M_1 = M_f(\Lambda_0)$ is $f(x) = x^3 (x-1)^2$.

Theorem~\ref{multiplicity:theorem:obstruction-3} is an instance of the following observation: 
Let $M = \bmat \mu_{i,j} \emat$ be an $m\times (n+1)$ multiplicity matrix 
such that $\colsum_j(M) = \sum_{i=1}^m \mu_{i,j} = n-j$ for some $j \in \{0,1,\ldots, n-1 \}$.   
If $M = M_f(\Lambda)$ for some monic polynomial $f(x)$ of degree $n$, then 
\[
f^{(j)}(x) = (n)_j \prod_{i=1}^{m} \left( x - \lambda_i \right)^{\mu_{i,j}} 
\]
where $\Lambda = \left( \lambda_1, \lambda_2, \ldots, \lambda_m\right)$.  
This constrains the shape of the multiplicity matrix
for $f(x)$ with respect to $\Lambda$, but, for $|\Lambda| = m \geq 3$,  not uniquely, 
even if $j = 0$. 
Here are some examples.  

\bex 
For $m = 4$ and $n=8$,  the $4 \times 9$ multiplicity matrices 
\[
N_1 = \bmat 
1 & 0 & 0 & 0 & 0 & 0 & 0& 0 & 0 \\
1 & 0 & 0 & 0 & 0 & 0 & 0& 0 & 0 \\
1 & 0 & 0 & 0 & 0 & 0 & 0& 0 & 0 \\
5 & 4 & 3 & 2 & 1  & 0 & 0& 0 & 0 
\emat
\]
and
\[
N_2 = \bmat 
1 & 0 & 1 & 0 & 0 & 0 & 0& 0 & 0 \\
1 & 0 & 0 & 0 & 0 & 0 & 0& 0 & 0 \\
1 & 0 & 0 & 0 & 0 & 0 & 0& 0 & 0 \\
5 & 4 & 3 & 2 & 1  & 0 & 0& 0 & 0 
\emat 
\]
satisfy 
\[
\colsum_0(N_1) = \colsum_0(N_2) = 8. 
\]
\eex

We have 
\[
N_1 = N_{f_1}(\Lambda_1)  
\]
for  
\[
\Lambda_1 = (0,1,2,3)
\]
and
\begin{align*}
f_1(x) & = x(x-1)(x-2)(x-3)^5 \\ 
& = x^8 -18x^7+ 137x^6 -570 x^5 + 1395x^4 -1998 x^3  +1539 x^2 - 486 x.  
\end{align*}
We have 
\[
N_2 = N_{f_2}(\Lambda_2). 
\]
for 
\[
\Lambda_2 = (0,2,3,-6)
\]
and
\begin{align*}
f_2(x) & = x(x-2)(x-3)(x+6)^5 \\ 
& = x^8 + 25x^7+216x^6 + 540 x^5 - 2160 x^4 -11664 x^3 + 46656 x.
\end{align*}

\bex  
 There are seven equivalence classes of $4 \times 5$ multiplicity matrices 
 $M = \bmat \mu_{i,j} \emat$ with $\mu_{i,0} =   1$ for $i \in \{1,2,3,4 \}$.  These matrices are 
\[
P_1 = \bmat
1 & 0 & 0 & 0 & 0 \\
1 & 0 & 0 & 0 & 0 \\
1 & 0 & 0 & 0 & 0 \\
1 & 0 & 0 & 0 & 0 
\emat,
\qquad 
P_2 = \bmat
1 & 0 & 0 & 1 & 0 \\
1 & 0 & 0 & 0 & 0 \\
1 & 0 & 0 & 0 & 0 \\
1 & 0 & 0 & 0 & 0 \emat,
\]
\[
P_3 = \bmat
1 & 0 & 1 & 0 & 0 \\
1 & 0 & 0 & 0 & 0 \\
1 & 0 & 0 & 0 & 0 \\
1 & 0 & 0 & 0 & 0 
\emat, 
\qquad 
P_4 = \bmat
1 & 0 & 2 & 1 & 0 \\
1 & 0 & 0 & 0 & 0 \\
1 & 0 & 0 & 0 & 0 \\
1 & 0 & 0 & 0 & 0 
\emat,
\]
\[
P_5 = \bmat
1 & 0 & 1 & 0 & 0 \\
1 & 0 & 0 & 1 & 0 \\
1 & 0 & 0 & 0 & 0 \\
1 & 0 & 0 & 0 & 0 
\emat, 
\qquad 
P_6 = \bmat
1 & 0 & 1 & 0 & 0 \\
1 & 0 & 1 & 0 & 0 \\
1 & 0 & 0 & 0 & 0 \\
1 & 0 & 0 & 0 & 0 
\emat,
\]
\[
P_7 = \bmat
1 & 0 & 1 & 0 & 0 \\
1 & 0 & 1 & 0 & 0 \\
1 & 0 & 0 & 1 & 0 \\
1 & 0 & 0 & 0 & 0 
\emat. 
\]
For  $k \in \{1,2,3,4,5\}$,  if $f(x)$ is a monic quartic polynomial 
such that   $P_k=M_f(\Lambda)$ with respect to $\Lambda = (0, a, b, c)$,  
then 
\begin{align*}
f(x) & = x(x-a)(x-b)(x-c).
\end{align*} 
\eex

It is straightforward to check that  
\begin{align*}
P_1 = M_f(\Lambda_1) & \qquad \text{for $\Lambda_1 = ( 0,1,2,3)$} \\
P_2  = M_f(\Lambda_2) &\qquad \text{for $\Lambda_2 = ( 0,1,2,-3)$} \\
P_3  = M_f(\Lambda_3) & \qquad \text{for $\Lambda_3 = ( 0,-3,4,12)$},
\end{align*}

Let  $f(x)$ be a monic quartic polynomial such that $P_4 = M_f(\Lambda)$ 
for some sequence $\Lambda$.  
By Theorem~\ref{multiplicity:theorem:affine-01}, we can assume that 
$\Lambda = (0,1, \rho_1, \rho_2)$.   
The multiplicity vector for $f(x)$ at $x=0$ is
$\bmat  1 & 0 & 2 & 1 & 0 \emat$ and so 
\[
f(x) = x^4 - ax = x\left(x^3-a\right)
\]
for some nonzero element $a$.  
Because $f(1) = 0$, it follows that $a=1 $  
and so   
\[
f(x) = x^4 - x = x(x-1)(x^2 +x+1). 
\] 
If  $x^2+x+1=0$ has no solution in $K$, then the matrix $P_4$ is not 
the multiplicity matrix of a polynomial in $K[x]$.  
If the polynomial $x^2+x+1$ splits in $K$ with zeros $\rho_1 = (-1+\sqrt{-3})/2$ 
and $\rho_2 = (-1- \sqrt{-3})/2$, 
then $M = M_f(\Lambda)$ for $\Lambda = ( 0,1,\rho_1,\rho_2)$.

Let  $f(x)$ be a monic quartic polynomial such that $P_5 = M_f(\Lambda)$ 
for the sequence $\Lambda = (0,1, \alpha_1, \alpha_2)$.  
The multiplicity vector for $f(x)$ at $x=0$ is
$\bmat  1 & 0 & 1 & 0 & 0 \emat$ and so 
\[
f(x) = x^4 +ax^3 + bx 
\]
for some nonzero elements $a$ and $b$.  
The multiplicity vector for $f(x)$ at $x=1$ is
$\bmat  1 & 0 & 0 & 1 & 0 \emat$ and so 
\[
f(x) = (x-1)^4 + c(x-1)^2 + d(x-1) 
\]
for some nonzero elements $c$ and $d$.  
We have 
\[
f(x) = x^4 +ax^3 + bx =  (x-1)^4 + c(x-1)^2 + d(x-1) 
\]
if and only if 
\begin{align*}
f(x) & = x^4 -4x^3 + 3x \\ 
& =  (x-1)^4 -6(x-1)^2  -5(x-1) \\
& = x(x-1)(x^2-3x-3).
\end{align*}
If  $x^2-3x-3=0$ has no solution in $K$, then the matrix $P_4$ is not 
the multiplicity matrix of a polynomial in $K[x]$.  
If the polynomial $x^2-3x-3$ splits in $K$ with zeros $\alpha_1 = (3 + \sqrt{21})/2$ 
and $\alpha_2 = (3 - \sqrt{21})/2$, 
then $M = M_f(\Lambda)$ for $\Lambda = ( 0,1,\alpha_1,\alpha_2)$. 
 
 Consider the $2 \times 5$ multiplicity matrix 
 \[
 M = \bmat 
1 & 0 & 1 & 0 & 0 \\
1 & 0 & 1 & 0 & 0 \\
\emat. 
 \]
Let $\Lambda_0= (0,1)$.  The polynomial 
\begin{align*}
 f(x) & = x^4  -2x^3+ x \\ 
 & = (x-1)^4  +2(x-1)^3 -(x-1) \\
 & = x(x-1)(x^2 - x - 1) 
\end{align*}
 is the unique monic quartic polynomial such that $M = M_f(\Lambda_0)$. 
 The zeros of this polynomial are $0,1, (1+\sqrt{5})/2$ and $(1- \sqrt{5})/2$. 
If the field $K$ does not contain the elements 
$\beta_1 = (1+\sqrt{5})/2$ and $\beta_2 = (1- \sqrt{5})/2$, 
then the multiplicity matrices $P_6$ and $P_7$ are not multiplicity matrices 
of polynomials in $K[x]$.  
If the field $K$ does contain $(1+\sqrt{5})/2$ and $(1- \sqrt{5})/2$, 
then the multiplicity matrix for the polynomial $ f(x) = x^4  -2x^3+ x$ 
with respect to $\Lambda = (0, 1, \beta_1, \beta_2)$ is $P_6$.  
It follows that the multiplicity matrix $P_7$ 
is not the multiplicity matrix of a polynomial in $K[x]$. 

These examples indicate the dependence of the multiplicity matrix on the field $K$.

\section{Open problems}          \label{multiplicity:section:OpenProblems}

\benum
\item      \label{Multiplicity:problem:1}        
Let $M$ be an $m\times (n+1)$ multiplicity matrix.   
Let $\Lambda = ( \lambda_1, \lambda_2, \ldots, \lambda_m)$ 
be a sequence of elements of $K$.     
Does there exist a polynomial $f(x)$ of degree $n$ 
such that $M_f(\Lambda) = M$? 

\item                     \label{Multiplicity:problem:2} 
Let $M$ be an $m\times (n+1)$ multiplicity matrix.      
Do there exist a sequence  
$\Lambda = ( \lambda_1, \lambda_2, \ldots, \lambda_m)$ of elements of $K$  
and a polynomial $f(x) \in K[x]$ of degree $n$ 
such that $M_f(\Lambda) = M$? 

\item                   \label{Multiplicity:problem:3} 
Are there additional constraints on the coordinates of a multiplicity matrix 
that guarantee that the matrix is the multiplicity matrix of a polynomial 
with respect to some sequence $\Lambda$?

\item           \label{Multiplicity:problem:Np1} 
Given a sequence $\Lambda = ( \lambda_1, \lambda_2, \ldots, \lambda_m)$ 
of elements of $K$ and an $m\times (n+1)$ multiplicity matrix $M$, 
problem~\eqref{Multiplicity:problem:1} asks if there is a polynomial $f(x)$ 
such that $M = M_f(\Lambda)$. 
There is the following more general question:  Does there exist a polynomial $F(x)$ 
of degree $N = n+p$ for some nonnegative integer $p$ 
such that $\lambda_i$ is a zero of $F^{(j)}(x)$ 
of multiplicity $\mu_{i,j}$ for all $i \in \{1,2, \ldots, m\}$ and $j \in \{0,1,\ldots, n\}$?  
If so, the multiplicity matrix $M_F(\Lambda)$ of the polynomial 
$F(x)$ with respect to $\Lambda$ is an $m \times (N+1)$-matrix that is an extension 
of the matrix $M$ in the sense that $M$ is the  submatrix of $M_F(\Lambda)$ 
obtained by deleting the $p$ columns $n+1,\ldots, N$ of $M_F(\Lambda)$.  
This is the \emph{extension problem}.
An example was given in Section~\ref{multiplicity:section:intro}. 

\item                       \label{Multiplicity:problem:Np2} 
Equivalently, the extension problem asks if every $m\times (n+1)$ multiplicity matrix 
can be extended 
to an $(m,n+p)$ multiplicity matrix $M'$ that is the multiplicity matrix of a polynomial 
$F(x)$ with respect to the  sequence $\Lambda$.
If such an extension is possible, what is the smallest degree 
of a polynomial $F(x)$ that extends $M$?

\item
Let $K = \R$ or \C.  Is there an analytic function $F(x)$ that extends $M$  
in the sense that $\lambda_i$ is a zero of $F^{(j)}(x)$ 
of multiplicity $\mu_{i,j}$ for all $i \in \{1,2,\ldots, m\}$ and $j \in \{0,1,\ldots, n\}$? 

\item
The results in this paper are valid for polynomials with coefficients in all fields of 
characteristic 0.  It would be of interest to investigate multiplicity vectors 
and  multiplicity matrices in fields of nonzero characteristic $p$.

\item
Classify the forbidden interactions between zeros of polynomials.

\eenum

\textbf{ Acknowledgements.}  I thank Sergei Konyagin, Noah Kravitz, Moshe Newman, 
Kevin O'Bryant, and other participants in the New York Number Theory Seminar 
for useful remarks.  
I also thank the referee for a very careful reading of this paper.

\end{document}